\magnification=\magstep0
\input amstex
\documentstyle{amsppt}
\NoBlackBoxes
\NoRunningHeads
\nologo
\topmatter
\title{A convolution formula for the Tutte polynomial}
\endtitle
\author W. Kook, V. Reiner, and D. Stanton
\endauthor
\address School of Mathematics, University of Minnesota,
Minneapolis, MN 55455.
\endaddress
\email
kook\@math.umn.edu, reiner\@math.umn.edu, stanton\@math.umn.edu
\endemail
\endtopmatter
Let $M$ be a finite matroid with rank function $r$.
We will write $A \subseteq M$ when
we mean that $A$ is a subset of the ground set of $M$, and
write $M|_A$ and $M/A$ for the matroids obtained by restricting $M$ to $A$, 
and contracting $M$ on $A$ respectively. 
Let $M^*$ denote the dual matroid to $M$. (See \cite{1} for definitions).
The main theorem is

\proclaim{Theorem 1} The Tutte polynomial $T_{M}(x,y)$ satisfies
$$
T_M(x,y)=\sum_{A \subseteq M} T_{M|_A}(0,y) T_{M/A}(x,0).
\tag1
$$
\endproclaim
First we define a convolution product and note a useful lemma.

Let $\Bbb{M}$ be the set of all isomorphism classes of finite matroids, 
and let $K$ be a commutative ring with 1. 
For any functions $f,g:\Bbb{M}\rightarrow K$, 
define $f\circ g:\Bbb{M}\rightarrow K$ by
$$
(f\circ g)(M)=\sum_{A\subseteq M} f(M|_A) g(M/A).
\tag2
$$ 
The convolution $\circ$ is associative, with identity element $\delta$,
$$
\delta(M)=
\cases
1 {\text{ if $M=\varnothing$,}}\\
0 {\text{ otherwise.}}
\endcases
$$  
Following Crapo \cite{2}, let $\zeta(x,y)(M)=x^{r(M)} y^{r(M^*)}$, where 
$K={\Bbb{Z}}[x,y]$.
\proclaim{Lemma 1} $\zeta(x,y)^{-1}=\zeta(-x,-y)$.
\endproclaim
\demo{Proof} Note that
$$
\align
(\zeta(x,y)\circ \zeta(-x,-y))(M)&
=\sum_{A \subseteq M} x^{r(M|_A)} 
y^{r((M|_A)^*)}(-x)^{r(M/A)}(-y)^{r((M/A)^*)}\\
&=x^{r(M)}y^{r(M^*)}\sum_{A \subseteq M} (-1)^{|M|-|A|}\\
&=\delta(M).
\endalign
$$
\enddemo
\demo{Proof of Theorem 1} The Tutte polynomial may be defined by \cite{1,2}
$$
T_{M}(x+1,y+1)=(\zeta(1,y)\circ \zeta(x,1))(M),
\tag3
$$
so also
$$
\align
T_{M}(x+1,0)&=(\zeta(1,-1) \circ \zeta(x,1))(M),\\
T_{M}(0,y+1)&=(\zeta(1,y) \circ \zeta(-1,1))(M). 
\endalign
$$
Therefore
$$
\align
\sum_{A \subseteq M} T_{M|_A}(0,y+1) T_{M/A}(x+1,0)&=(\zeta(1,y)\circ 
\zeta(-1,1))
\circ (\zeta(1,-1)\circ \zeta(x,1))(M) \\
&=\zeta(1,y)\circ (\zeta(-1,1)
\circ \zeta(1,-1))\circ \zeta(x,1)(M)\\
&=\zeta(1,y)\circ \zeta(x,1)(M)\\
&=T_{M}(x+1,y+1),
\endalign
$$
where the third equality is by Lemma 1.
\qed
\enddemo

\noindent
{\bf{Remark 1.}}
 
Note that Theorem 1 can be rewritten as
$$
T_M(x,y)=\sum_{{\text{isthmus-free flats $V$}}} 
T_{V}(0,y) T_{M/V}(x,0).
\tag4
$$
This is because when $A \subseteq M$ is not a flat,
$M/A$ contains a loop $e$ and
$$
T_{M/A}(x,0)=\left[ y \, T_{(M/A)-e}(x,y) \right]_{y=0}=0.
$$
Similarly if $A$ contains an isthmus $e$, then 
$$
T_{M|_A}(0,y)=\left[ x \, T_{(M|_A)/e}(x,y) \right]_{x=0}=0.
$$

\noindent
{\bf{Remark 2.}}
 
Theorem 1 can also be proven using Tutte's original
definition of the Tutte polynomial involving
{\it basis activities} \cite{1,2}: for any ordering of the ground set of $M$, 
$$
T_M(x,y):=\sum_{\text{bases }B\text{ of }M} x^{|IA_M(B)|} y^{|EA_M(B)|}
$$
where here $IA_M(B)$ (resp. $EA_M(B))$ denotes the set of internally
(resp. externally) active elements of $M$ with respect to the
base $B$.

\newpage
\noindent
Theorem 1 in \cite{3}  asserts that any base $B$ can be 
uniquely decomposed $B=B_1 \cup B_2$ with $B_1 \cap B_2 =\varnothing$ and
$$
IA_V(B_1)=EA_{M/V}(B_2) = \varnothing
$$
where $V$ is the flat $\overline{B_1}$ spanned by $B_1$.
It turns out that in this decomposition one furthermore has
$$
IA_M(B)=IA_{M/V}(B_2), \qquad EA_M(B)=EA_{V}(B_1). 
\tag5
$$
We omit the details of this verification, which are straightforward.
Given this, one then has
$$
\align
T_M(x,y)
&=\sum_{\text{bases }B\text{ of }M} x^{|IA_M(B)|} y^{|EA_M(B)|}\\
&=\sum_{\text{flats }V\text{ of }M} 
  \sum_{\text{bases }B_1\text{ of }V \atop
        \text{ with }IA_V(B_1)=\varnothing}
  \sum_{\text{bases }B_2\text{ of }M/V \atop
        \text{ with }EA_{M/V}(B_2)=\varnothing}
x^{|IA_{M/V}(B_2)|} y^{|EA_V(B_1)|} \\
&=\sum_{\text{flats }V\text{ of }M} 
  \left(\sum_{\text{bases }B_1\text{ of }V \atop
        \text{ with }IA_V(B_1)=\varnothing} y^{|EA_V(B_1)|}  \right)
  \left(\sum_{\text{bases }B_2\text{ of }M/V \atop
        \text{ with }EA_{M/V}(B_2)=\varnothing} x^{|IA_{M/V}(B_2)|} \right)
\\
&=\sum_{\text{flats }V\text{ of }M} T_V(0,y) \, T_{M/V}(x,0).
\endalign
$$

\noindent
{\bf{Remark 3.}}
  
The version (4) of Theorem 1 may also be proven by deletion-contraction,
as we now explain.  Recall \cite{1} that the Tutte polynomial
is characterized by the following three properties.
\roster
\item"(i)" $T_M(x,y) = x$ if $M$ consists of a single isthmus,
and $T_M(x,y)=y$ if $M$ consists of single loop.
\item"(ii)" $T_{M_1 \oplus M_2}(x,y) = T_{M_1}(x,y) \cdot T_{M_2}(x,y).$
\item"(iii)" $T_M(x,y) = T_{M-e}(x,y) + T_{M/e}(x,y)$ if
$e$ is neither an isthmus nor a loop of $M$.
\endroster

Let $T'_M(x,y)$ be the right side of (4), and we must show that it
also satisfies (i),(ii),(iii).  Properties (i),(ii) are straightforward
and omitted.  To show (iii), fix an element $e$ which is neither an
isthmus nor a loop of $M$, and then use property (iii)
for $T_M(x,0), T_M(0,y)$ to write
$$
\align
T'_M(x,y)&=\sum_{{\text{isthmus-free flats $V$}}} 
              T_{V}(0,y) T_{M/V}(x,0) \\
         &=\sum_{{\text{i.f. flats $V$} \atop e \in V}} 
              T_{V}(0,y) T_{M/V}(x,0)
         +\sum_{{\text{i.f. flats $V$} \atop e \not\in V}}        
              T_{V}(0,y) T_{M/V}(x,0)\\
         &=\sum_{{\text{i.f. flats $V$} \atop e \in V}} 
              T_{V-e}(0,y) T_{M/V}(x,0)
         +\sum_{{\text{i.f. flats $V$} \atop e \in V}}        
              T_{V/e}(0,y) T_{M/V}(x,0)\\
         &+\sum_{{\text{i.f. flats $V$} \atop e \not\in V}} 
              T_{V}(0,y) T_{M/V-e}(x,0)
         +\sum_{{\text{i.f. flats $V$} \atop e \not\in V}} 
              T_{V}(0,y) T_{(M/V)/e}(x,0)\\
         &=\sum_{V, V-e \text{ both i.f.} \atop e \in V} 
              T_{V-e}(0,y) T_{M/V}(x,0)\\
         &+\sum_{{\text{i.f. flats $V$} \atop e \in V}} 
              T_{V/e}(0,y) T_{M/V}(x,0)\\
         &+\sum_{{\text{i.f. flats $V$} \atop e \not\in V}} 
              T_{V}(0,y) T_{M/V-e}(x,0)\\
         &+\sum_{\text{i.f. flats $V$} \atop 
                 V\cup\{e\}\text{ a flat}} 
              T_{V}(0,y) T_{(M/V)/e}(x,0)
\tag6
\endalign
$$
where the last equality comes from the fact that $T_{V-e}(0,y)=0$
unless $V-e$ is isthmus-free, and dually $T_{(M/V)/e}(x,0)=0$
unless $V\cup\{e\}$ is a flat of $M$.

On the other hand, we wish to show that the above sum is
the same as
$$
\align
&T'_{M-e}(x,y) + T'_{M/e}(x,y)\\
   &=\sum_{\text{i.f. flats }W \atop \text{ of }M-e}
        T_{(M-e)|_W}(0,y) T_{(M-e)/W}(x,0)
    +\sum_{\text{i.f. flats }W \atop \text{ of }M/e}
        T_{(M/e)|_W}(0,y) T_{(M/e)/W}(x,0)\\
   &=\sum_{\text{i.f. flats }W\text{ of }M-e, \atop
           W\text{ not a flat of }M}
        T_{(M-e)|_W}(0,y) T_{(M-e)/W}(x,0)\\
    &+\sum_{\text{i.f. flats }W\text{ of }M/e, \atop
           W\text{ not a flat of }M}
        T_{(M/e)|_W}(0,y) T_{(M/e)/W}(x,0)\\
    &+\sum_{\text{i.f. flats }W\text{ of }M-e, \atop
           W\text{ a flat of }M}
        T_{(M-e)|_W}(0,y) T_{(M-e)/W}(x,0)\\
    &+\sum_{\text{i.f. flats }W\text{ of }M/e, \atop
           W\text{ a flat of }M}
        T_{(M/e)|_W}(0,y) T_{(M/e)/W}(x,0)
\tag7
\endalign
$$
The terms $W$ in the sums on the right-hand side of equation (7)
biject with the terms $V$ in the sums on the right-hand side of
equation (6) as follows:  in the first sum $W=V-e$, in the
second sum $W=V/e$, in the third sum $W=V$ and in the fourth
sum $W=V$.  We leave it to the reader to check that this gives
a bijection of the terms which shows the equality of the
right-hand sides in (6) and (7).  The only tricky point here is
in the fourth sum, where one must note that not only are $W,V$
equal as subsets of the ground sets of $M/e, M$ respectively, but
also the flats $W, V$ of $M/e, M$ are isomorphic as matroids,
due to the fact that $e$ is an isthmus of $V\cup\{e\}$.

\vskip .1in
\noindent
{\bf{Remark 4.}}
 
Lemma 1 can be used to prove other convolution identities.
For example, if we define
$$
\rho(x,y,z,w)(M):= (\zeta(z,y) \circ \zeta(x,w))(M)
$$
then equation (3) implies 
$$
\align
T_M(x,y) &= \rho(x-1,y-1,1,1)(M) \\
T_M(0,y) &= \rho(-1,y-1,1,1)(M) \\
T_M(x,0) &= \rho(x-1,-1,1,1)(M)
\endalign
$$
and Theorem 1 is the specialization $z=w=1$ of the more general
identity
$$
\align
\rho(x-1,y-1,z,w) 
&= \zeta(z,y-1) \circ \zeta(x-1,w)\\
&= \zeta(z,y-1) \circ \zeta(-1,1) \circ \zeta(1,-1) \circ \zeta(x-1,w)\\
&= \rho(-1,y-1,z,1) \circ \rho(x-1,-1,1,w)
\endalign
$$

As another example, of the use of Lemma 1, one can
start with equation (2) and multiply both sides by
$\zeta(-1,-y)$.  Using the notation $T(x,y)(M):=T_M(x,y)$, we obtain 
$$
\zeta(-1,-y) \circ T(x+1,y+1) =\zeta(x,1)
$$
which gives an apparently new recursion for the Tutte polynomial
$$
T_M(x,y) =  
 (x-1)^{r(M)} - \sum_{\varnothing \neq A \subseteq M} (-1)^{r(M|_A)}
                   (1-y)^{r(M|_A^*)} T_{M/A}(x,y).
$$

\noindent
{\bf{Remark 5.}}
 
The convolution product defined by equation (2) 
suggests a certain coalgebra (actually a Hopf algebra) 
naturally associated with matroids, first considered in \cite{4, \S 15}. 
Let $A$ be a free $K$-module with basis ${\Bbb M}$ equal
to the isomorphism classes of finite matroids $[M]$. The coproduct 
$\Delta:A\rightarrow A\otimes A$ is defined $K$-linearly by
$$
\Delta([M])=\sum_{A \subseteq M} [M|_A]\otimes [M/A],
$$
and the product $\mu:A\otimes A\rightarrow A$ is defined 
$K$-linearly by
$$
\mu([M]\otimes [M'])=[M\oplus M'].
$$
Define a bigrading on $A$ by setting the bidegree of $[M]$ to be
$(r(M),r(M^*))$. One can check that this makes
$A$ a co-associative, commutative, bigraded, connected, 
Hopf algebra over $K$, whose unit $\eta:K\rightarrow A$ is
$\eta(1)=[\varnothing]$, and whose 
co-unit $\epsilon:A\rightarrow K$ is $\epsilon([M])=\delta_{M,\varnothing}$.
If $\phi:A\rightarrow A$ is the involution $\phi([M])=[M^*]$ extended
$K$-linearly to all of $A$, then one can check that the identity
$M^*|_{M-A} \cong (M/A)^*$ leads to the equation 
$$
\Delta\circ \phi=(\phi\otimes \phi)\circ \Delta^{op}.
$$
Therefore $\phi^*:A^*\rightarrow A^*$ is an algebra 
anti-automorphism.  Note that $\phi$ also exchanges the bigrading
in the sense that if $a$ has bidegree $(s,t)$ then $\phi(a)$ has
bidegree $(t,s)$.  

Motivated by this, let $A$ be any co-associative, bigraded, connected
coalgebra over $K$ with coproduct $\Delta$ and co-unit $\eta$,
having a distinguished $K$-basis of bihomogeneous elements ${\Bbb M}$.
Let $\circ$ denote the product dual to $\Delta$ in the dual algebra $A^*$, 
and $\phi:A \rightarrow A$ be any involution which exchanges
the bigrading and such that $\phi^*:A^*\rightarrow A^*$
is an anti-automorphism.
Define $\zeta \in A^*$ by $\zeta(x,y)(M) = x^s y^t$ for all
$M \in {\Bbb M}$ having bidegree $(s,t)$.  We can then
define a {\it Tutte functional} $T(x,y) \in A^*$
by $T(x,y) = \zeta(1,y-1) \circ \zeta(x-1,1)$.
One can then check that the
familiar Tutte polynomial identity \cite{1}
$$
T_{M^*}(x,y) = T_M(y,x)
$$
has the counterpart 
$$
\phi^*(T(x,y))=T(y,x)
$$
which follows formally from the assumed properties of $\phi$.

Furthermore, the proof of Lemma 1 actually shows the following in
this context: 
$$
\text{If }\zeta(1,1)^{-1} = \zeta(-1,-1)
\text{ then }\zeta(x,y)^{-1} = \zeta(-x,-y)^{-1}.
$$
Consequently, if we impose the extra condition on $A$ that
$\zeta(1,1)^{-1} = \zeta(-1,-1)$, then the counterpart
to Theorem 1
$$
T(x,y) = T(0,y) \circ T(x,0)
$$
ensues as a formal consequence.


\noindent
{\bf Acknowledgments.}
The authors thank Richard Stanley for suggesting the search for
a Tutte polynomial analogue to equation (2.1) of \cite{3}.


\Refs

\ref
\no 1
\by T. Brylawski and J. G. Oxley
\paper The Tutte polynomial and its applications
\inbook Matroid Applications
\ed N. White
\publ Cambridge Univ. Press
\publaddr Cambridge
\yr 1992
\endref

\ref
\no 2
\by H. Crapo
\paper The Tutte polynomial
\jour Aequationes Mathematicae
\vol 3
\yr 1969
\pages 211-229
\endref

\ref
\no 3
\by W. Kook, V. Reiner, and D. Stanton
\paper Combinatorial Laplacians of matroid complexes
\paperinfo preprint
\yr 1997
\endref


\ref
\no 4
\by W. Schmitt
\paper Incidence Hopf Algebras
\yr 1994
\jour J. Pure Appl. Algebra
\vol 96
\pages 299-330
\endref

\endRefs
\enddocument
\end